\documentstyle[12pt]{article}

 \def\bbb{\bf} \def\blb{\small \bf} 
 \def\varnothing{\emptyset} 

\def\medskipamount{12pt}

\newcommand{\bit}[1]{\section{#1}\setcounter{equation}{0} }
\newcommand{\bitt}[1]{\subsection{#1} }

\renewcommand{\theequation}{\thesubsection .\arabic{equation}}

\newcommand{\re}[1]{{\bf (\ref{#1})}}

\catcode`\@=\active
\catcode`\@=11
\def\@eqnnum{\hbox to .01pt{}\rlap{\bf \hskip -\displaywidth(\theequation)}}
\catcode`\@=12

\newenvironment{s}[1]
{ \addvspace{\medskipamount} \refstepcounter{equation}
\noindent {\bf (\theequation) #1.} \begin{em}}
{\end{em} \par \addvspace{\medskipamount} }

\newenvironment{r}[1]
{ \addvspace{\medskipamount} \refstepcounter{equation}
\noindent {\bf (\theequation) #1.} }
{\par \addvspace{\medskipamount} }

\begin{document}

\catcode`\@=\active
\catcode`\@=11
\newcommand{\nc}{\newcommand}


\nc{\vars}[2]
{{\mathchoice{\mbox{#1}}{\mbox{#1}}{\mbox{#2}}{\mbox{#2}}}}
\nc{\Aff}{\vars{\bbb A}{\blb A}}
\nc{\C}{\vars{\bbb C}{\blb C}}
\nc{\G}{\vars{\bbb G}{\blb G}}
\nc{\Hyp}{\vars{\bbb H}{\blb H}}
\nc{\N}{\vars{\bbb N}{\blb N}}
\nc{\Pj}{\vars{\bbb P}{\blb P}}
\nc{\Q}{\vars{\bbb Q}{\blb Q}}
\nc{\R}{\vars{\bbb R}{\blb R}}
\nc{\V}{\vars{\bbb V}{\blb V}}
\nc{\Z}{\vars{\bbb Z}{\blb Z}}


\nc{\oper}[1]{\mathop{\mathchoice{\rm #1}{\rm #1}
{\scriptstyle \rm #1}{\scriptstyle \rm #1}}\nolimits}
\nc{\Aut}{\oper{Aut}}
\nc{\Def}{\oper{Def}}
\nc{\End}{\oper{End}}
\nc{\Hilb}{\oper{Hilb}}
\nc{\Hom}{\oper{Hom}}
\nc{\diag}{\oper{diag}}
\nc{\Fl}{\oper{Fl}}
\nc{\Gr}{\oper{Gr}}
\nc{\NS}{\oper{NS}}
\nc{\Par}{\oper{Par}}
\nc{\Pic}{\oper{Pic}}
\nc{\Proj}{\oper{Proj}}
\nc{\Quot}{\oper{Quot}}
\nc{\Spec}{\oper{Spec}}


\nc{\GL}[1]{{\rm GL(#1)}}
\nc{\PSL}[1]{{\rm PSL(#1)}}
\nc{\PGL}[1]{{\rm PGL(#1)}}
\nc{\SL}[1]{{\rm SL(#1)}}

\def\commrect#1#2#3#4#5#6#7#8{%
\begin{center}%
\begin{picture}(130,90)%
\put(120,70){\vector( 0,-1){50}}%
\put(10,80){\vector( 1, 0){100}}%
\put(0,70){\vector( 0,-1){50}}%
\put(10,10){\vector( 1, 0){100}}%
\put(115,80){\makebox(0,0)[l]{#2}}%
\put(5,80){\makebox(0,0)[r]{#1}}%
\put(115,10){\makebox(0,0)[l]{#4}}%
\put(5,10){\makebox(0,0)[r]{#3}}%
\put(-3,50){\makebox(0,0)[r]{#5}}
\put(123,50){\makebox(0,0)[l]{#6}}
\put(60,3){\makebox(0,0)[c]{#8}}
\put(60,88){\makebox(0,0)[c]{#7}}
\end{picture}\end{center}}

\def\commtriang#1#2#3#4#5#6{%
\begin{center}\begin{picture}(120,80)%
\put(55,70){\vector(-1,-2){30}}
\put(65,70){\vector(1,-2){30}}
\put(30,5){\vector(1,0){60}}
\put(60,75){\makebox(0,0)[c]{#1}}
\put(25,5){\makebox(0,0)[r]{#2}}
\put(95,5){\makebox(0,0)[l]{#3}}
\put(60,0){\makebox(0,0)[c]{#6}}
\put(37,43){\makebox(0,0)[r]{#4}}
\put(83,43){\makebox(0,0)[l]{#5}}
\end{picture}\end{center}}


\nc{\down}{\Big\downarrow}
\nc{\beqas}{\begin{eqnarray*}}
\nc{\beqa}{\begin{eqnarray}}
\nc{\beq}{\begin{equation}}
\nc{\bl}{\vskip 1.2ex }
\nc{\eeqas}{\end{eqnarray*}}
\nc{\eeqa}{\end{eqnarray}}
\nc{\eeq}{\end{equation}}
\nc{\fp}{\mbox{     $\Box$}}
\nc{\half}{\frac{\scriptstyle 1}{\scriptstyle 2}}
\nc{\m}{{\bf m}}
\nc{\mod}{/ \! \! /}
\nc{\pf}{{\em Proof}}
\nc{\sans}{\backslash}
\nc{\st}{\, | \,}
\nc{\cC}{{\cal C}}
\nc{\hcC}{{\hat {\cal C}}}

\catcode`\@=12

\newpage
\title {The Topology of Conjugate Varieties}
\author {David Reed\\
Mathematical Institute\\
24 - 29 St Giles'\\
Oxford OX1 3LB\\
UK\\}
\date{October 1994}
\maketitle

\begin{abstract}  Serre \cite{Se:64}  and Abelson \cite{Ab:74}  have produced
examples of conjugate algebraic
 varieties which are not homeomorphic.  We show that if the field of definition
of a polarized projective variety coincides with its field of moduli then all
of its conjugates have the same topological type.  This immediately extends the
class of varietie
s known to possess invariant topological type to all canonically embedded
varieties.  We also show that (normal) complete intersections in projective
space and, more generally  in homogeneous varieties, satisfy the condition.
\end{abstract}

\bit{Introduction}
If $V$ is an algebraic variety defined over $k$, a finitely generated extension
of $\Q$, for
 each embedding $\sigma :k \rightarrow \C$ we can extend scalars to form the
complex
algebraic variety $V_{\sigma}$ defined by the following cartesian
square (base change or extension of scalars)
\[
\begin{array}{ccc}
V_{\sigma} := V \times _{\Q} \C & \rightarrow & V\\
\downarrow &    & \downarrow \\
\Spec \C & \stackrel{\sigma}{\rightarrow} & \Spec k
\end{array}
\]
The complex points of this variety $V_{\sigma}(\C)$ form a topological space
and the
topological type of two such spaces, $V_{\sigma}(\C)$ and $V_{\tau}(\C)$ for
two
 different embeddings $k \rightarrow \C$ can  be compared.  We will refer to
varieties
 $V_{\sigma}$  and $V_{\tau}$ obtained in this manner as {\em conjugate
varieties}.

Serre \cite{Se:64} and Abelson \cite{Ab:74} have produced examples of conjugate
varieties whose complex points constitute non-homeomorphic topological spaces.
Since the publication of these papers there has been little published work in
this area.  The principal result of the research reported upon here is a
sufficient condition for
the topological spaces of complex points of conjugate varieties to be
homeomorphic.  We begin by recalling definitions due to
Matsusaka-Shimura-Koizumi \cite{Ko:72}.

\begin{r}{Definition}
For a divisor $X$ on a projective variety $V$ define the class ${\cal P}(X)$ to
be the class of all divisors $X'$ on $V$ such that there are integers $m, n$
with $mX \equiv nX'$ (algebraic equivalence).  If ${\cal P}(X)$ contains an
ample divisor then it
is called a {\em polarization} (this term is also applied to the ample divisor
in the class).  An isomorphism of projective varieties $f:V \rightarrow W$ is
said to be a isomorphism of polarized projective varieties if there are
polarizations ${\cal P}$,
${\cal P}'$ on $V$, and $W$ respectively such that  the map on divisors induced
by $f$ takes ${\cal P}$ to ${\cal P}'$
 \end{r}

\begin{r}{Definition}
The field of moduli for  polarized projective variety $(V,{\cal P})$ is the
field $k$ such that for $\sigma \in \Aut(\C)$, $\sigma$ is in fact in $\Aut(\C
/k)$  if and only if $V^{\sigma} \simeq V$ as polarized varieties.\end{r}

Compare this field to the {\em field of definition} of the variety which is the
field $K$ such that $\sigma \in \Aut(\C)$ is in fact in $\Aut(\C /k)$ iff and
only if $V^{\sigma}=V$.

\bl
A discussion of when  fields of moduli exist in general for varieties can be
found in \cite{Ko:72}.  An example of a variety (a hyperelliptic curve) whose
field of moduli differs from its field of definition can be found in
\cite{Sh:72}.

\begin{s}{Theorem}
\label{big'}
If $V$ is a polarized projective variety defined over $k$, a finitely generated
extension of $\Q$, and if the field of moduli for $V$ coincides with $k$, then
the topological type of $V_{\sigma}(\C)$ is independent of $\sigma$.\end{s}

The proof of this Theorem relies on a strengthened form of Thom's
 stratified Isotopy Theorem which is given in \S 2 below.  The Theorem itself
is proved in \S3.

It was previously known that certain types of varieties whose topology is
rather easily described, such as non-singular curves, abelian varieties, $K-3$
surfaces and simply connected non-singular surfaces of general type, have
topological types which do n
ot vary under conjugation of their fields of definition.  The above Theorem
extends this list to include all canonically embedded varieties.

In \S 4, the criterion given above will be used to show that (normal) complete
intersections in projective space, and, more generally, in homogeneous
varieties, also belong on this list by showing that they satisfy the condition
of the Theorem as well. In
 \S 5 we point out a few of the many questions that remain open in this area.

\bit{Stratified Isotopy Theorems and the Topology of Conjugate Varieties}
\bitt{Preliminaries from Algebraic Geometry}

Our intention is to review the proof of the ``stratified isotopy Theorem'' as
it applies to algebraic varieties with a view to establishing that the
stratification described by the Theorem can be defined without extension of the
base fields of the varieti
es involved.   This material is essentially contained in \cite{Ve:76}.

In the course of the  analysis we will frequently rely on two sets of
well-known and indeed basic facts from algebraic geometry which are stated here
with an emphasis on the relevant fields of definition.
The first set of facts deals with the singular locus of a variety defined over
a field $k$ of characteristic 0 (here we do not need any restriction on the
nature of the extension $k/\Q$).

 Let $X$ be a variety of dimension $n$ over a characteristic zero field $k$,
then for any point $x \in X$ the following are equivalent:
\begin{enumerate}
\item  $\Omega^1_{X,x}$ (the module of differentials at $x$) is a free module
of rank $n$ over the local ring ${\cal O}_{X,x}$ of $X$ at $x$;

\item  ${\cal O}_{X,x}$ is a regular local ring.

\end{enumerate}

and, if either of these conditions obtains at $x$ we say $X$ is {\em smooth} at
$x$.

\begin{s}{Fact}
\label{sm}
 There is an everywhere dense Zariski open set $U \subset X$ which is smooth
and the Zariski closed set $X-U$ is defined over $k$.
\end{s}

{\bf Remark}:  For a  discussion of fields of definition for arbitrary subsets
of schemes see [{\bf EGA IV}, \S 4.8].

The second set of facts is just Hironaka's well-known resolution of
singularities.  Once again our only restriction is that we work over
characteristic 0 fields.

Let $X$ be variety defined over $k$ and ${\cal J}$ a coherent sheaf of ideals
defining a closed sub-scheme $D$ then we make the usual:

\begin{r}{Definition}

A blow-up of $X$ at $D$, otherwise known as a monoidal transformation of $X$
with center $D$, is a pair $(P, f)$ consisting of a variety $P$ and a morphism
$f:P \rightarrow X$ such that $f^{-1}({\cal J})$ is an invertible sheaf on $P$
and, for any other p
air $(P',f')$ with $f':P' \rightarrow X$ and $f'^{-1}({\cal J})$ an invertible
sheaf on $P'$, there is a unique morphism $g:P' \rightarrow P$ such that
\commtriang{$P'$}{$P$}{$X$}{$g$}{$f$}{$f'$}

commutes.\end{r}

A general procedure for constructing $P$ is to define
\[
P:=\Proj(\oplus_{d=0}^{\infty} {\cal J}^d)
\]
where we set ${\cal J}^0={\cal O}_X$.  There is a natural map $P \rightarrow X$
(given by ${\cal O}_X \rightarrow \oplus {\cal J}^d$) and the universal
property is proved in \cite{H:77} Ch II, \S 7.  In particular the field of
definition for $P$ is just t
he field of definition of $D$ or, equivalently, of ${\cal J}$.
\bl
Hironaka has shown,

\begin{s}{Theorem}
\label{m1}
Let $X$ be a variety defined over $k$, characteristic 0, then there is a closed
subscheme $D$ of $X$ such that:
\begin{enumerate}
\item  the set of closed points of $D$ is the singular locus of $X$; and
\item  if $f:{\tilde X} \rightarrow X$ is the monoidal transform of $X$ at $D$
then $\tilde X$ is smooth.
\end{enumerate}
\end{s}

\pf.  \cite{Hi:64} (Main Theorem 1)  \fp

For our purposes we note in particular that since $D$ is defined over $k$ by
Fact \re{sm} above we have that $\tilde X$ and $f$ are defined over $k$ as
well.

\begin{r}{Definition}
A divisor with normal crossings $D$ in a smooth variety $X$ is a divisor such
that for any $x \in D \subset X$ with local ring ${\cal O}_{X,x}$ and maximal
ideal ${\bf m}_{X,x}=(z_1, \dots ,  z_n)$, each component of $D$ passing
through $x$ is described b
y precisely one ideal $(z_i)$.
\end{r}

\begin{s}{Theorem}
\label{ms}
Let $X$ be a smooth variety defined over $k$, $W$ a nowhere dense sub-scheme of
$X$, then there exists a finite set of monoidal transforms
\[
f_i:X_{i+1} \rightarrow X_i
\]
with smooth centers $D_i$, for $0 \leq i < r$ and $X_0=X$ such that
\begin{enumerate}
\item  $X_r$ is smooth;
\item if $\bar{f_i}$ is the composition of the $f_j$ for $0 \leq j < i$, then
$D_i \subset {\bar {f_i}}(W)$ for all $i$; and
\item  ${\bar{f_r}}^{-1}(W)$ is an invertible sheaf whose support is a divisor
with normal crossings.
\end{enumerate}
\end{s}

\pf.  \cite{Hi:64} (Cor 3, to Main Theorem II).  \fp

Since we know that the singular locus of a variety $X$ over $k$ is nowhere
dense and hence its inverse image under the monoidal transform $f$ from Theorem
\re{m1} is nowhere dense, we can summarize the above by saying that Hironaka's
resolution of singula
rities starts with an arbitrary variety $X$, defined over $k$ characteristic 0,
and produces a smooth variety $X'$ and a morphism $f:X' \rightarrow X$ such
that the inverse image of the singular locus of $X$ becomes a divisor with
normal crossings in $X'$
 and $X'$, $D$ and $f:X \rightarrow X$ are defined over $k$ as well.

\bitt{Stratifications and stratified isotopy in the Real Analytic Category}

      The most natural setting for the study of stratifications of singular
spaces is the category of real analytic subspaces of smooth (real analytic)
manifolds and proper maps between them.  A brief sketch of the aspects of the
theory used below is give
n here.  The application to complex algebraic varieties follows.  The best
current reference is \cite{G-M:88}.

Let $M$ be a real analytic manifold, $Z \subset M$ a closed subset and
\[
Z = \bigcup_{i\in S}S_i
\]
($S$ a partially ordered set) a decomposition of $Z$ as a union of a locally
finite collection of disjoint locally closed ``pieces'' or ``strata''
satisfying the {\em boundary} consition
\[
S_i \cap {\bar{S_j}} \neq \varnothing \Leftrightarrow S_i \subset {\bar
{S_j}}\Leftrightarrow i=j\;\;or\;\;i<j
\]
(in the last case we also write $S_i < S_j$).  Such a decomposition is called a
{\em Whitney Stratification} if and only if it also satisfies:
\begin{enumerate}
\item  each $S_i$ is smooth (not necessarily connected), and
\item each pair $(S_i,S_j)$ satisfies the {\em a} and {\em b} conditions,
namely, if we have a collection of points $\{x_i\} \subset S_i$ such that
$\{x_i\} \rightarrow y \in S_j$ and another set of points $\{y_i\} \subset S_j$
with $\{y_i\} \rightarrow y
$ such that the secant lines ${\overline{x_iy_i}}\rightarrow l$ and the tangent
planes $T_{x_i}S_i \rightarrow \tau$ then we have

 \begin{itemize}
\item {\em a}:  $T_yS_j \subset \tau$; and
\item {\em b}:  $l \subset \tau$
\end{itemize}
\end{enumerate}

These conditions ensure that the pieces $S_i$ ``fit together'' well at an
infinitesimal level (see \cite{B-C-R:87} for examples).  The conditions are
local and can be tested by taking local coordinates in $M$ about $y$.  The
validity of the conditions is
independent of the choice of coordinate system.  It is a theorem
(Hironaka-Hardt) that any subanalytic manifold admits such a stratification.

If a map behaves well with respect to stratifications we say it is a {\em
stratified} map.  Specifically, let $Y_1 \subset M_1, \;\; Y_2 \subset M_2$ be
Whitney stratified subsets of manifolds $M_1, \; M_2$, and let $f:M_1
\rightarrow M_2$ be a real analy
tic map such that $f\!\mid Y_1$ is proper and $f(Y_1) \subset Y_2$, then $f$ is
{\em stratified} if for each stratum $A\subset Y_2$ we have $f^{-1}(A)$ a union
of connected components of strata of $Y_1$, say $f^{-1}(A) = \cup S_i$ and $f$
takes each $S_i$
 submersively to $A$ (surjection on tangent spaces).  There are two key results
on stratified maps which are often referred to as the $1^{st}$ and $2^{nd}$
(stratified) Isotopy Theorems.

\begin{s}{Theorem}  For $Z \subset M$ a Whitney stratified subset of a real
analytic manifold, $f:Z \rightarrow \R^n$ proper and such that the restriction
to each stratum
\newline $f\!\mid A : A \rightarrow \R^n$ is a submersion, then there is a
stratum preserving homeomorphism $h:Z \rightarrow \R^n \times (f^{-1}(0) \cap
Z)$ such that

\commrect{$Z$}{$\R^n\times(f^{-1}(0)\cap
Z)$}{$\R^n$}{$\R^n$}{$f$}{$pr_1$}{$h$}{$id$}

\noindent commutes.  In particular, the fibers of $f\!\mid Z$ are homeomorphic
by a stratum preserving homeomorphism.
\end{s}

\begin{s}{Theorem}
$A \subset M$, $B \subset N$ subanalytic subsets of real analytic manifolds,
$F: A \rightarrow B$ a proper subanalytic map.  Then there exist
stratifications $S$, $T$ of $A$, $B$ into smooth subanalytic manifolds such
that $f$ is stratified with respect t
o $S$ and $T$.  Furthermore, given any locally finite collection of subanalytic
subsets $\cal C$ of $A$ (resp. ${\cal D}$ of $B$ we can choose $S$ (resp. $T$)
such that each elements of ${\cal C}$ (resp. ${\cal D}$) is a union of strata
of $S$ (resp. $T$)
{}.

\end{s}

By the $1^{st}$ isotopy Theorem one obtains local topological triviality of the
$f$ along connected components of strata of $B$.

For further discussion and guidance to the literature see \cite{G-M:88} Part I,
Chapter 1, pp. 36 - 44.

We  wish to employ this theory in the context of complex algebraic varieties
and to take our stratifications to be constructible sets whose fields of
definition we can control.  Following a suggestion of Bernstein, Beilinson,
Deligne \cite{BBD:81} Chapter
 6, we find that such a version of  stratified isotopy theory has been given by
Verdier.

\bitt{Whitney Stratifications \`{a} la Verdier}

We now define a notion of Whitney stratification which is adapted to algebraic
varieties.  The properties {\em a} and {\em b} above will be replaced by a
single property {\em w} which is also local in nature.  Hence we will always
assume that {\em smooth}
 complex algebraic varieties have been equipped with coordinate charts given by
their underlying real analytic manifold structure.  Nothing will depend on the
choice of coordinates (see comments below).

We use the following notion of distnace between sub-vector spaces in a finite
dimensional Euclidean space $E$, ${\delta}(F,G)$
 defined by
\[{\delta}(F,G) := \sup_{\small{\begin{array}{ccc}  x\! & \in & F \\ \parallel
x \parallel & = & 1 \end{array}}}\!dist (x,G)\]

In particular, $\delta(F,G)=0 \Rightarrow F \subset G$.

\begin{r}{Definition}

A Verdier-Whitney stratification of (the complex points of) a complex algebraic
variety $V$
where $V$ is a $k$-variety
of finite type,  $k$ a finitely generated extension of $\Q$, is a finite
disjoint partition
 of $V$ by smooth constructible sets $A_i$
\[V=\bigcup_{{i}=1}^n A_{i}\]
 such that
 \begin{enumerate}
\item the ``boundary property'' holds, namely ${\overline{A_{i}}} \cap A_{j}
\neq \varnothing$ implies ${\overline {A_{i}}} \supset A_{j}$ and
\item if ${\overline {A_{i}}} \supset A_{\j}$ with $i \neq j$ then the pair
$(A_{i},A_{j})$ satisfies
the following property ``{\em w}" at every point $y \in A_{i}$:

Consider $A_{\alpha}$ and $A_{\beta}$ as real analytic manifolds  and take
coordinate patches around $y$ to
some
Euclidean space $E$, then there exists a neighborhood $U \subset E$ of (the
image of)
$y$
and a positive real number $C$ such that $\forall x \in U \cap A_{i}$ and $y'
\in U \cap A_{j}$
(here $A_i$ and $A_j$ are taken to mean the images of some small open subsets
around $y$ in $E$) we have
\[{\delta}(T_{y'}{A_{j}},T_x{A_{i}}) \leq C \parallel x - y' \parallel\]
where $T_x{A_{\alpha}}$ is the tangent plane to $A_{\alpha}$ at $x$ and
$\delta$ is as defined above

\end{enumerate} \end{r}

A number of remarks are called for here.

As we are restricting ourselves to algebraic varieties, it is sufficient to
consider
 stratifications with finite collections of subsets.  The analytic cases
require infinite
collections of subsets.  This permits a certain amount of simplification in the
definition and the subsequent arguments.
It also permits us to speak of the (common) field of definition of the
stratification as being the smallest field containing the fields of definition
of the $A_{\alpha}$.

On the other hand the condition {\em w} (so-called by Verdier) replaces the
more familiar
 conditions `{\em a} and `{\em b} above. Condition {\em w} implies condition
{\em a} simply because it is a uniform version
 of it but {\em w} does not imply {\em b} in general (consider the logarithmic
spiral at 0).
The key fact however is that this implication does hold when $A_{\alpha}$ is a
smooth
 subanalytic subspace of a real analytic space and $A_{\beta}$ is a smooth
analytic
subspace of
$\overline {A_{\alpha}}$ (Kuo).

Verification of property {\em w} does not depend on the choice of coordinates
(for this it is important that the strata $A_i$ are required to be smooth).
The next key fact is that using resolution of singularities, we can stratify
arbitrary complex algeb
raic varieties in much the same way as real analytic manifolds.

The first Theorem we will require is,

\begin{s}{Theorem}
\label {3b}
If $V$ is a complex algebraic variety as above and $V_{\beta}$ a finite family
of
 constructible subsets of $V$ then there is a Verdier-Whitney stratification of
$V$ such
that each $V_{\beta}$ is obtained as the union of strata.   The stratification
is defined over the (common) field of definition of $V$ and the
$V_{\beta}$.\end{s}

The proof of this is based on

\begin{s}{Theorem}
\label {3a}
$V$ as above, $M$,  $M'$ smooth, connected, locally closed subsets
such that $M \cap M' = \varnothing$, $M' \subset {\overline M}$ and ${\overline
M} - M'$ is
 closed
(all for the Zariski topology), then there is a Zariski open $Y \subset M'$
containing all of the points $y \in M'$ such that $(M, M')$ has the property w
at $y$ and
$M' - V$ is Zariski closed.  $Y$ is defined over the (common) field of
definition of $V$, $M$ and $M'$. \end{s}

\pf\quad  of \re{3a}.  When $M$, $M'$ are locally closed smooth {\em
subanalytic subspaces of a second countable real analytic space} $X$ a
corresponding Theorem is proved by
\begin{enumerate}
\item defining a subset $V \subset M'$ by removing ``bad points'' to arrive at
an open subanalytic subset of $X$ which is dense in $M'$, and then
\item showing that $(M, M')$ has property {\em w} at all points of $V$ by
taking coordinate charts in which $M'$ is an open subset of an affine space
$F$. There is an affine space $G$ such that locally $F\oplus G=X$, there is a
``blow-up'' $W$ of $X$ in w
hich $M'$ is described as a divisor with normal crossings $\prod_1^qz^{n_i}_i$
and there is a map $\pi:W \rightarrow X=F\oplus G$.  Additional work involving
an analysis of the matrix representation of $d\pi$ then gives the desired
result.

In what follows we show how $V$ is defined in the complex algebraic case, show
that $V$ is a Zariski open and describe its field of definition.  We do not
review the proof that $(M, M')$ has the property {\em w} at all points of $V$
since we are principal
ly interested in showing that the stratifications can be taken to be algebraic
and have the right fields of definition.  The reference for all omitted parts
of proofs is \cite{Ve:76}.

We use
 resolution of singularities to find a smooth complex algebraic variety $W$ and
a proper
 morphism $\pi :W \rightarrow V$ with $\pi(W)={\overline M}$ and
$\pi^{-1}(M')=D \subset W$ is a divisor whose singularities are at worst normal
crossings.  By \re{m1} and \re{ms} resolution of singularities  takes place
without extension of the field of
definition so that both $W$ and $D$ are defined over the same field as $V$.
The next step is to produce the desired open set by removing ``bad points'', in
this case, the normal crossings singularities.

Consider the subset $D_q \subset D$ of points of $D$ where at least $q$
irreducible
local components of $D$ meet (this set is empty for $q \gg 0$) let ${\tilde
{D_q}}$ be the normalization of $D_q$ (separating the points lying on the
various components)
and let $i:{\tilde{D_q}}\rightarrow D_q$ be the normalization map.
${\tilde{D_q}}$ is
 a smooth algebraic variety defined over the same field as $D$ so we can
consider
\[d({\pi}\circ i_q):\; \Omega^1_{D_q} \rightarrow (\pi \circ
i_q)^*\Omega^1_{\overline M}\]
where $\Omega^1$ is the sheaf of differentials.  We now have one more
correction to make,
 namely we must consider the points where this map is not
 surjective so that we do not have a submersion.  These form a  Zariski closed
subset $R_q \subset D_q$ [EGA IV, \S. 17.15.13] defined by a Jacobian condition
and hence this singular locus is defined over the same field as $D_q$.

$i_q(R_q)$ is closed since $i_q$ is finite and hence a closed map (the
``going-down" Theorem) and the collection of
$i_q(R_q)$ is finite so that $S:=\pi(\bigcup_q i_q(R_q)) \subset {\overline M}$
is Zariski
 closed.  We thus have that \[Y:=M' \cap (M - S) \subset {\overline M}\] is a
Zariski open
 subset of $M$ defined without extension of the base field such that $Y$ is
dense in $M'$
 and it is smooth by construction.  The proof  of \re{3a} now proceeds by
considering the local analysis of
 the smooth real analytic varieties undrlying $M$, $M'$ and $Y$ as very briefly
described above.\fp

\pf\quad of Theorem \re{3b}.  We need two Lemmas.

\begin{s}{Lemma}
\label{311}
Let $X$ be an algebraic variety, $Y_{\beta}$ a finite family of constructible
subsets of $X$, then there is another finite family of subsets of $X$,
$B_{\alpha}$
 satisfying: \begin{enumerate}
\item for all $\alpha$, ${\overline B}_{\alpha}$ and ${\overline
B}_{\alpha}-B_{\alpha}$ are
Zariski closed and the $B_{\alpha}$ are connected and smooth;
\item the $B_{\alpha}$ partition $X$ and each $Y_{\beta}$ is the union of a
collection
of $B_{\alpha}$;
\item ${\overline B}_{\alpha} \cap B_{\beta} \neq \varnothing \Rightarrow
B_{\beta} \subset {\overline B}_{\alpha}$
\end{enumerate} and the $B_{\alpha}$ can be defined without extending the
(common)
 field of definition of the $Y_{\beta}$. \end{s}

\pf.   Since the $Y_{\beta}$ are locally closed we have that the sets
${\overline Y}_{\beta}$ and ${\overline Y}_{\beta} - Y_{\beta}$ are Zariski
closed and we replace the family $Y_{\beta}$ with the family of Zariski closed
sets $\{X,{\overline Y}_{\beta
},{\overline Y}_{\beta}-Y_{\beta}\}$ which we continue to refer to as
$Y_{\beta}$.

Let $\cal F$ be the largest collection of Zariski closed subsets of $S$ which
is such that: \begin{itemize}
\item for all $\beta$, $Y_{\beta} \in {\cal F}$;
\item ${\cal F}$ is closed under intersection;
\item for all $Z \in {\cal F}$ the irreducible components of $Z$ are in ${\cal
F}$; and
\item $Z \in {\cal F}$ implies that the set of singular points of $Z$, $sing\;
Z \in {\cal F}$.
\end{itemize}

We can get at least one such collection with these properties by taking the
family of sets
 consisting of the $Y_{\beta}$ and their singular loci $(Y_{\beta})_{sing}$ and
then
 closing this collection under taking of irreducible components (defined over
the field of
definition of the $Y_{\beta}$) and intersections.  In particular $\cal F$ can
be
constructed without extending the field of definition of the $Y_{\beta}$.

Then we define
\[B_{\alpha}=Z_{\alpha} - \bigcup_{\stackrel{Z_{\beta}\subset
Z_{\alpha}}{Z_{\beta} \neq Z_{\alpha}}}Z_{\beta}\]
where the $Z_{\gamma}$ are the irreducible sets in $\cal F$.  It is clear that
the
 $B_{\alpha}$ have the properties set out in the Lemma.\fp

\begin{s}{Lemma}
\label{3l2}
Let $V \subset X$ be a constructible set which is connected and smooth and let
$Z
\subset {\overline V} - V$ be Zariski closed, then there is a (finite)
partition
$B_{\alpha}$ of $Z$ such that such that for all $\alpha$, ${\overline
B}_{\alpha}$,
${\overline B}_{\alpha} - B_{\alpha}$ are Zariski closed in $X$, the
$B_{\alpha}$ are
smooth and connected and the pairs $(V, B_{\alpha})$ have property {\em w}.
The
 $B_{\alpha}$ are defined over the (common) field of definition of $X$, $V$ and
$Z$
\end{s}

\pf.  By induction.  The statement is true for $Z= \varnothing$.  Assume $Z
\neq \varnothing$ so we apply Lemma \re{311} to get a finite collection of
smooth connected
 subsets $U_{\alpha} \subset Z$, $U_{\alpha} \cap U_{\beta} = \varnothing$
defined
over the field of definition of $Z$such that the ${\overline U}_{\alpha}$,
${\overline U}_{\alpha} -U_{\alpha}$, $Z-U_{\alpha}$ are all Zariski closed and
 $Z-\bigcup_{\alpha}U_{\alpha}$ is a Zariski closed set of lower dimension.  By
Theorem \re{3a} there is an open dense $W_{\alpha} \subset U_{\alpha}$ such
that the
 ${\overline U}_{\alpha} - W_{\alpha}$ are Zariski closed and the pairs
 $(V, W_{\alpha})$ have the property {\em w}.

Now $Z_1:=Z - \cup W_{\alpha}$ is Zariski closed with dimension lower than $Z$
so
 we can apply the induction hypothesis to it.  Combine the $B_{\alpha}$ thus
obtained
with the $W_{\alpha}$ to get a new collection of $B_{\alpha}$.  Since the
 $W_{\alpha}$ provided by Theorem \re{3a} are constructed without extending
fields of
definition we are done.\fp

Returning to the Proof of Theorem \re{3b} we start by restricting to the case
of $V$
irreducible and proceed once again by induction.

For $V= \varnothing$ the Theorem is trivially true so suppose $V \neq
\varnothing$ and
replace the $V_{\beta}$ in the statement of the Theorem with a the family $\{V,
{\overline V}_{\beta}, V - V_{\beta}, {\overline V}_{\beta} - V_{\beta}\}$ as
in the
 proof of Lemma \re{311} so that this new family (which we still call
$Y_{\beta}$) is
 made up of Zariski closed sets.

Since $V$ is irreducible there is a ${\beta}_0$ such that $Y_{{\beta}_0}$ is
open and
dense in $V$.  Hence $V_1:=V-Y_{{\beta}_0}$ is closed and is of lower dimension
than
$V$.  We also clearly have that $\beta \neq {\beta}_0 \Rightarrow Y_{\beta}
\subset
 X_1$.

Let $B_{\alpha}$ be a partition of $V_1$ coming from Lemma \re{3l2} and apply
 Lemma \re{311} to the $Y_{\beta}$ and $B_{\alpha}$ together.  This produces a
common refinement $\{C_{\gamma}\}$ which still has property {\em w} because, in
general, if $M$, $M'$ are locally closed subsets of an algebraic variety $X$
which are
 smooth with $M' \subset \overline M$ and $M' \cap M = \varnothing$, then if
there is a
 locally closed and smooth $M'' \subset M'$, the pair $(M, M'')$ will have the
property
{\em w} if the pair $(M,M')$ does.

Now assume for the moment the following \newline
{\bf Claim}:  If $V = V_{\alpha}$ (finite union) and the $V_{\alpha}$ are
Zariski closed,
then if Theorem \re{3b} is true for the $V_{\alpha}$ it is true for $V$.

$V_1$ is a finite union of irreducibles of dimension lower than $V$ so apply
the
induction hypothesis to $V_1$ and the $C_{\gamma}$ to get a Whitney
stratification of
$V_1$ and add $Y_{{\beta_0}}$

If we can now prove the {\bf Claim} we have just made we will both complete the
proof
for the irreducible case and for the general case as well.  So let $V=Y \cup Z$
be a union
of irreducibles.  Apply the result in the irreducible case to $Z$, $Y_{\beta}
\cap Z$, $Y
\cap Z$ to get a Whitney stratification of $Z$ such that $Y_{\beta} \cap Z$ and
$Y \cap Z$ are unions of strata.
Now apply the result again to $Y$, $Y_{\beta} \cap Y$ and the $A_{\alpha}$ such
that
 $A_{\alpha} \subset Y \cap Z$ to get a Whitney stratification of $Y$ such that
the
$A_{\alpha}$ and $Y_{\beta} \cap Y$ are unions of strata.  Take $B_{\beta}$ and
those $A_{\alpha}$ such that $A_{\alpha} \subset Z - (Y \cap Z)$.  Note that
property
{\em w} still holds by the remark made above. Continue by induction. Nothing in
any of
these procedures requires an extension of fields of definition.\fp

\bitt{Stratified Morphisms}

The next ingredient is the demonstration that quite general algebraic-geometric
morphisms behave well with respect to
 stratifications.

Recall that a morphism of stratified spaces $f:X \rightarrow Y$ is called a
{\em stratified} morphism if it is proper and if the inverse image of a stratum
of $Y$ under $f$ is a union of strata of $X$ and each component of these strata
is mapped sunmersiv
ely (as a real analytic manifold) to $Y$.

We now have a version of the $2^{nd}$ isotopy theorem for complex algebraic
varieties.

\begin{s}{Theorem}(Verdier)
\label{3c}
If $f:X \rightarrow Y$ is a morphism of complex algebraic varieties and $f$ is
proper
then there are Verdier-Whitney stratifications $S$ and $T$ of $X$ and $Y$,
defined over the (common) field of
definition of $X$,$Y$ and $f$, such that $f$ is transverse to $S$ and $T$
\end{s}

\pf.

{\em Step 1}  We first show that if $X \rightarrow Y$ is proper and $S$ a
stratification of $X$ such as given in Theorem \re{3b}, then there is a Zariski
open $U \subset f(X) \subset Y$ which is smooth in $Y$ and (Zariski) dense in
$f(x)$ such that $f\!\m
id f^{-1}(U) \rightarrow U$ takes the connected components of$S \cap f^{-1}(U)$
submersively to $U$.

Set
\[X_q=\bigcup_{dim\:S_{\alpha} \leq q}S_{\alpha}\]
this is closed, smooth, constructible subset of $X$ and we write $f_q$ for
$f\mid X_q:X_q \rightarrow Y$.

Once again, let $R_q$ be the set of points in $X_q$ where
\[df_q: \Omega^1_{X_q} \rightarrow f^*_q(\Omega^1_Y)\]
is not onto and set
\[U_q=Y-Y_{reg}\cap (\bigcup_{q'\leq q}f_{q'}(R_{q'}))\]
where $Y_{reg}$ is the set of regular points of $Y$ (as noted above this set is
defined
 without extension of the field of definition in characteristic 0).  The proof
now
 proceeds by induction on $q$ to show that
\[f_q\mid U_q:f^{-1}_q(U_q) \rightarrow U_q\]
maps the connected components of to $S\cap f^{-1}_q(U_q)$ to $f_q(U_q)$ where
we always have ${U_q}$ Zariski dense in $Y$.  Note that $Y-U_q$
is automatically a Zariski closed subset since the $f_{q'}(R_{q'})$ are.

So let $S_{\alpha} \subset X_{q}$ be a stratum.  If  dim $S_{\alpha} < q$ then
there is a $q'<q$ such that $f_{q'}$ maps the connected components of $S \cap
f^{-1}_{q'}(U_{q'})$ to $U_{q'}$.  Thus we can assume dim$S_{\alpha}=q$ so
$S_{\alpha} \subset X_{
q}$ is
 Zariski open and smooth and is a connected component of $W=X_q - X_{q-1}$
which
is
also open and smooth.

 Now $R_q \cap f^{-1}(U_{q-1}) \subset W$ since
$\forall q' < q$
\[coker\:(df_{q'}) \rightarrow coker\:(df_q)\]
is surjective over points of  $X_{q'}$  as $f_q$ agrees with $f_{q'}$ on
$X_{q'}$.
 $f_q(R_q) \cap U_{q-1}$ is a Zariski closed subset and $U_{q-1} - (f_q(R_q)
\cap U_{q-1})$ is dense in $U_{q-1}$  and
$f_q \mid U_q$ takes $S \cap f^{-1}_q(U_q)$ submersively to $f_q(U_q)$.

 But $f:X \rightarrow Y$ is
 proper so $f_q(R_q)$ is empty for $ q \gg 0$ and we can form $U:= \bigcap_q
U_q$.  This is dense in $Y$ with
$Y - U$ an algebraic Zariski closed set and by construction $f$ takes the
connected components of $S\cap f^{-1}(U)$ submersively to $U$.

{\em Step 2}  Now take $f:X \rightarrow Y$ with stratifications $S$ and $T$
respectively as guaranteed by Theorem \re{3b}.  By step 1 we find a Zariski
open $U \subset f(X)$ smooth in $Y$ and dense in $f(X)$ with $f\!\mid
f^{-1}(U)$ submersive on connecte
d components.  We now re-stratify $Y$ per Theorem \re{3b} with $U$ a union of
strata.  Take $f^{-1}$ of these strata and restratify $X$.  $f$ now behaves as
we want on $f^{-1}(U)$.  Consider $f\!\mid X-f^{-1}(U)$, find a $U' \subset
f(X-f^{-1}(U))=f(X)-f(
U) \subset Y$ which is Zariski open, smooth in $Y$ and dense in $f(X)-f(U)$ and
repeat the above restratification process.  Since $dim\: U' < dim \:U$ the
procedure terminates after a finite number of steps.  Using Theorem \re{3b} we
will always have that
 the $U$, $U'$,  $\dots$ will be unions of strata and similarly for the
$f^{-1}(U)$, $f^{-1}(U')$,  $\dots$ and $f$ will clearly take each connected
component of the stratifications of the $f^{-1}(U)$, $f^{-1}(U')$,  $\dots$ to
$f(U)$, $f(U')$, $\dots$

 \fp

\bitt{Stratified Isotopy}

Lastly we state a version of the $1^{st}$ isotopy Theorem:

\begin{s}{Theorem}(Thom)
\label{Th}
$X$, $Y$ are real analytic spaces, $S$ and $T$  stratifications,
$f:X \rightarrow Y$ proper and submersive on the connected components of the
strata of $X$.
Set $y_0 \in Y$.  Write  $X_0=f^{-1}(y_0)$,  $S_0=X_0 \cap S$.
  Then there is an open neighborhood (in the complex topology) $y_0 \in V
\subset Y$ and a homeomorphism
$\phi:(f^{-1}(V),\;S \cap f^{-1}(V)) \rightarrow (X_0 \times V,\;S_0 \times V)$
preserving the stratifications and compatible with projections to $V$.\end{s}

\pf.  Classically this is proved using techniques from differential topology
and is quite difficult.  Using the condition {\em w} in place of the more
standard {\em a} and {\em b} conditions Verdier is able to give a fairly
self-contained proof in a few p
ages \cite{Ve:76}.  Nonetheless we will pass this over in silence since our
objective is not to see how the result is proved but rather to show how it can
be applied to give useful transversality properties for the algebraically
defined stratifications de
scribed above.  See \cite{Ve:76} for the missing details. \fp

When combined with Verdier's results this gives:

\begin{s}{Corollary}
\label{313}
  Let $X \stackrel{f}{\rightarrow} Y$ be a proper morphism of algebraic
varieties,  then the topological type of the fibers of $f$ over a connected
component of a stratum of $T$ is constant. \end{s}

\pf.
 By Theorem \re{3c} there are stratifications $S$ and $T$ of $X$ and $Y$
respectively, defined over the common fields of definition of $X$, $Y$ and $f$,
 such that the inverse image of any stratum of $T$ is a union of strata of $S$
$f^{-1}(T_{\alpha}) =
\cup S_{\beta}$ and each connected component of a stratum of $S$ is mapped
submersively onto a stratum of $T$.  Thus only the last statement requires
discussion.  Consider a connected component of a stratum of $T$, and call it
$W$. Partition $W$ into subs
ets such that the topological type of the fibers of $f$ are constant on each
member of the partition.  The sets partitioning $W$ are then open by Theorem
\re{Th} and they are disjoint by construction.  Since $W$ is connected only one
of the sets in the pa
rtition is non-empty.  \fp

\bit{Principal Results}

\bitt {A Sufficient Condition for Topological Stability Under Conjugation}

The results described above can be applied to give a sufficient criterion for
an algebraic
 variety and its conjugates to have the same topological type.

\begin{s}{Theorem}
\label{main}
Let $V$ be a $k$-variety, $k$ a finitely generated extension of $\Q$ and
suppose there exists a family  $f:Y \rightarrow B$, that is, a proper morphism
of complex algebraic varieties such that:
\begin{enumerate}
\item all of the conjugate complex algebraic varieties $V_{\sigma}$ are
isomorphic as $k$-varieties to fibers of $f$
(in other words for each $\sigma : k \rightarrow \C$ there is a point
$b_{\sigma} \in B$ and $k$-isomorphism $V_{\sigma} \simeq f^{-1}(b_{\sigma})$);
and
\item   $f: Y \rightarrow B$ arises by base extension from $\Q$,  that is,
there are $\Q$ varieties $Y_{/{\Q}}$ and $B_{/{\Q}}$ and  a morphism
$f_{/{\Q}}: Y_{/{\Q}} \rightarrow B_{/{\Q}}$ such that
\commrect{$Y \simeq Y_{/{\Q}} \times \C$}{$Y_{/{\Q}}$}{$B \simeq B_{/{\Q}}
\times \C$}{$B_{/{\Q}}$}{$f$}{$f_{/{\Q}}$}{${\beta}_Y$}{${\beta}_B$}
commutes,
\end{enumerate}

then the topological type of $V_{\sigma}(\C)$ is independent of $\sigma$.
\end{s}

\pf.  By Corollary \re{313} we can stratify $Y$ and $B$ with stratifications
$S$ and $T$
defined over $\Q$ so that the map $f$ is topologically locally trivial over
each connected component of the strata.  We need only show therefore that the
points $b_{\sigma}$,$b_{\tau}$ corresponding
to conjugate varieties $V_{\sigma}$, $V_{\tau}$ must lie in a single connected
 component of a stratum of $S$.

Since $v_{\sigma}:=f^{-1}(b_{\sigma})$, $v_{\tau}:=f^{-1}(b_{\tau})$ are
$k$-isomorphic to varieties which differ only by conjugation of their field of
definition these fibers of $f$  are mapped to the same subscheme $v_{/{\Q}}$ of
 $Y_{/{\Q}}$ by ${\beta}_Y$.  This is in turn
mapped to a subscheme of $B_{/{\Q}}$; call it $b_{/{\Q}}$. Now we claim
$v_{/{\Q}}$ and hence $b_{/{\Q}}$ are irreducible.  If not, the $v_{\sigma}$
divide up into subsets which are interchanged by some $\phi \in \Aut(\C)$ (each
$v_{\sigma}$ is irreducibl
e since it is isomorphic to a variety $V_{\sigma}$. Thus the largest
irreducible closed subschemes of $Y_{\Q}$ containing the images of these
subsets are not equal to $Y_{\Q}$.  But $Y_{\Q}$ is irreducible since $Y$ is,
so we have a contradiction.  Thus $
b_{/{\Q}}$ is
defined by a sheaf of prime ideals ${\cal P}$ with local ring ${\cal O}_{\cal
P}$ and
residue field ${\bf k}_{\cal P}$.

Take the inverse image of $b_{/{\Q}}$ under ${\beta}_B$ in $B \simeq B_{/{\Q}}
\times \C$.  Call this $b_{/{\C}}$.  The points
$b_{\sigma}$ and $b_{\tau}$ lie in $b_{/{\C}}$ by
commutativity of the diagram in Theorem \re{main} and call the residue fields
of these points ${\bf k}_{\sigma}$ and ${\bf k}_{\tau}$.  Now ${\beta}_B$ is an
open map  so that we have ${\bf k}_{\sigma} \simeq {\bf k}_{\cal P} \otimes \C$
and
similarly for ${\bf k}_{\tau}$.

The stratification $T$ of $B$ is defined over $\Q$ and we claim that points of
a subscheme of $B$, defined over $\Q$ and with residue fields isomorphic to
$k_{\cal P}$ over $\C$ must lie in a single irreducible component of a stratum
of $T$.

Assume not.  By Verdier's results $b_{/{\C}}(\C)$, the set of complex points of
 $b_{/{\C}}$ is a union of strata $ \cup T_i$ and each $T_i$ is defined over
$\Q$.  Assume that we have $b_{\sigma} \in T_1$ and $b_{\tau} \in T_2$.  We may
assume $T_1 \cup T_2$ is all of $b_{/{\C}}(\C)$.  By the definition  of a
stratification, if ${\overline {T_1}} \cap T_2 \neq \varnothing$ then
${\overline {T_1}} \supset T_2$.  This implies that one of the points
$b_{\sigma}$ or $b_{\tau}$ lies in a zariski closed subset of $b_{\C}$ and
hence does not have residue field isomorphic to $k_{\cal P}$ over $\C$.  So we
must have ${\overline {T_1}} \cap
 T_2 = \varnothing$ and similarly ${\overline {T_2}} \cap T_1 = \varnothing$.
 But then
$b_{/{\Q}}$ is not irreducible. Thus the points $b_{\sigma}$ and $b_{\tau}$ lie
in a single irreducible component of the stratification of $B$. Finally, over
$\C$, irreducibility in the zariski topology implies connectedness in the
complex topology so  $b
_{\sigma}$ and $b_{\tau}$ lie in a single connected component of the
stratification and hence the topological types of $f^{-1}(b_{\sigma})$ and
$f^{-1}(b_{\tau})$ are the same.\fp

{\em Remark}:  Shimura has used the existence of non-homeomorphic conjugate
varieties to show that the irreducible components of Chow varieties are not
necessarily defined over $\Q$ \cite{Sh:68}.  I am grateful to Professor J-P.
Serre for providing this r
eference.

\bitt{Corollaries}

  Following, for example \cite{BBD:81}, it is easy to see that any complex
projective variety  $V$ (say defined over a field $k$) can be embedded in a
family
 defined over $\Q$.  One merely considers the coefficients $c_{\alpha \beta}$
of the
 homogeneous ideal defining $V$ as indeterminates.  This produces a family
$f: Y \rightarrow S$ over an affine base $S$ with the original variety $V$
isomorphic to
the fiber of $f$ over the point of $S$ corresponding to the $c_{\alpha \beta}$.
This family clearly contains all of the conjugates of $V$, since these are
obtained by conjugating the coefficients in its homogeneous ideal.
Furthermore,  since $S$ is affi
ne
this family is actually ``arises via base extension from  $\Q$" in the sense of
Theorem \re{main}.

Fortunately we cannot use this technique to show that all varieties have
invariant  topological type under conjugation due to the fact that the families
we obtain in this way may not be irreducible.  Consider as an example the
hyperelliptic curve ${\cal C
}$ constructed by Shimura \cite{Sh:72}
\[
y^2=a_0x^m + \sum^m_{r=1}(a_rx^{m+r} + (-1)^r{\bar {a_r}}x^{m-r})
\]

If we treat the $a_r$ and $\bar{a_r}$ as independent indeterminates $a_r$,
$b_r$ we get a family of curves $F$ and $\cal C$ lies over a point $p_0$ in the
locus where $a_r=\bar{b_r}$.  The family has two components which are
interchanged by complex conjug
ation. We cannot use it in Theorem \re{main} therefore to show that the
conjugates of $\cal C$ are homeomorphic (although, of course, this can be shown
in other ways).

If the field of moduli of $V$ coincides with its field of definition however,
we can use the BBD type family in Theorem \re{main} by virtue of,

\begin{s}{Proposition}
Let $V$ be a projective variety whose field of definition $k$, a finitely
generated extension of $\Q$, coincides with its field of definition.  Define
the family $f: Y \rightarrow S$ as above, by letting the coefficients
$\{c_{\alpha \beta}\}$ of the homo
geneous ideal of $V$ vary.  Then all of the conjugates of $V$ are isomorphic to
 fibers $f^{-1}(b_i)$ where the $b_i$ lie in a single connected component of a
Verdier-Whitney stratification of $S$ and hence have the same topological type
.\end{s}

\pf. Examining the proof of Theorem \re{main} we see that we did not need there
the full assumption that $Y$ is irreducible, but merely that the conjugates of
$V$ do not lie in separate irreducible components of $Y$.

Consider the action of $\sigma \in \Aut(\C)$ on $f:Y \rightarrow S$.  Since the
field of moduli of $V$ coincides with its field of definition $k$ we have that
$Y \rightarrow S$ is stable under $\sigma$ if and only if $\sigma \in \Aut (\C
/k)$ (that is, th
e field of definition of $f:Y \rightarrow S$ is $k$).  Let $Y=\cup Y_i$ be a
decomposition of $Y$ into irreducible components with fields of definition
$k_i$. If $Y_i \neq Y$ then $k \subset k_i$ and $k \neq k_i$.  Let the fibers
of $f$ corresponding to t
he conjugates $V_{\sigma}$ be $v_{\sigma}$ and suppose, for example, that
$v_{\sigma}$ and $v_{\tau}$ were in $Y_i$ and $Y_j$ respectively. Then there is
some $\phi \in \Aut(\C/ k)$ taking $Y_i$ to $Y_j$.  but $\phi$ fixes
$v_{\sigma}$ and $v_{\tau}$ so b
oth must be in $Y_i \cap Y_j$.  This argument applies to all of the
$v_{\sigma}$ and all components $Y_i$ hence all of the conjugates of $V$ can be
identified with fibers of $f$ lying in a single component $Y_i$ (we may pick
any component).  We now replac
e $Y$ with $Y_i$, which is irreducible, and apply the theorem.

\fp

\begin{s}{Corollary}
The topological type of a canonically embedded variety is invariant under
conjugation.\end{s}

\pf.  The canonical embedding is given over the field of definition and hence
the field of definition coincides with the field of moduli. \fp

This enlarges the class of varieties previously known to have conjugate
invariant topological type.  Further examples are given in \S 4 where we show
that complete intersections in homogeneous varieties have fields of moduli
which coincide with the fields
 of definition.

For a proof that Serre's original examples of non-homeomorphic conjugate
varieties do not have fields of definition which coincide with their fields of
moduli see \cite{Re:94}

\bit{Complete Intersection Type Varieties}

\bitt{Generalities}

 We can use Theorem \re{main} to exhibit further classes of
varieties which have topological type invariant under conjugation.  These
classes of varieties, which include (normal) complete intersection varieties in
projective
space,  may be parametrized by vector spaces of sections of vector
bundles and this linear structure provides a natural method of
descending from $\C$ to $\Q$.

The Deformation Theory of these varieties has been studied by a
number of authors including \cite{K-S:58}, \cite{S:75}, \cite{B:83} and
\cite{W:84}
in a series of papers with results extending from smooth hypersurfaces through
to the more general cases.  The most general result along these lines is:

\begin{s}{Theorem} (Wehler)
\label{weh}
Let $Z=G/H$ be a non-singular homogeneous complex variety, quotient of a
simple,
simply connected Lie group $G$ by a parabolic subgroup $H$,
$E=\oplus_{j=1}^r{\cal O}_Z(d_j)$ a vector bundle, $s \in H^0(Z,E)$ a
section and $X$ a complex variety described by the zero locus of $s$
such that codim$X=r$ and $X$ is not a $K$-3 surface.  Then the vector
space $H^0(Z,E)$ parametrizes a complete set of small deformations of
$X$ and these deformations are given by the family
\[ Y :=\{(z,s) \mid z\in Z, s\in H^0(Z,E), s(z)=0\} \rightarrow
H^0(Z,E)
\]
\end{s}

  If we are to apply the Theorem to this case we must show that the family $Y
\rightarrow H^0(Z, E)$ satisfies the conditions of Theorem \re{main}.  To do
this we reprove the theorem using algebraic geometric techniques to obtain a
family over $\Q$ which
will have the required properties.
The proofs given here therefore are
modelled  on Wehler's but adapted to algebraic rather
than analytic geometry.  The key in both the analytic and algebraic
approach is to establish a close connection between the deformation
theory of the objects and their Hilbert schemes as will be explained shortly.
For another approach to establishing the algebraic deformation theory of
complete intersections see \cite{Ma:68}.

\bitt{Comparison of Hilbert Scheme and Deformation Functors}

 Let $Z$ be an arbitrary  non-singular projective variety over a field
$k$ and $E \rightarrow Z$ an algebraic $k$ vector bundle over $Z$. Consider the
scheme $X
\subset Z $ defined by a global section $s_0 \in H^0(Z,E)$.  We construct this
scheme as
follows:  the section $s$
defines a map of
sheaves ${\cal O}_Z \rightarrow E$ sending the section ${\bf 1}$ of
${\cal O}_Z$ to the
section $s$ (we are abusing notation here by not distinguishing between the
vector bundle and its associated locally free sheaf).  There is a dual
map
$\check s:\check E \rightarrow {\cal O}_Z$  and $X$ is said to be
defined by $s$ if ${\cal O}_X$ fits into an
exact sequence of sheaves
\[\check E \stackrel{\check s}{\rightarrow} {\cal O}_Z \rightarrow
{\cal O}_X \rightarrow 0 \]
$X$ is sometimes referred to as the zero scheme of $s$.

 \begin{r}{Definition} The Hilbert functor of $X$ (in $Z$), is
the functor $\Hilb$ from $\cal C$, the category of local artin rings with
residue field $k$ to the category ${\em Sets}$ which assigns to each
object $A$ in $\cal C$ the set of schemes $Y$ which fit into the
following diagram

\commrect{$Z \supset X\simeq Y \times_k A$}{$Y \subset Z \times \Spec
A$}{$\Spec k$}{$\Spec
A$}{}{}{}{}
with $Y$ flat over $\Spec A$.\end{r}

\begin{r}{Definition}  The (affine) projective cone (hereinafter simply
``the cone'') $C_X$ of (or on) a projective
variety $X \subset \Pj^n$ is given by
\[
C_X:= \Spec ({\cal O}_X \oplus
{\cal O}_X(1) \oplus {\cal O}_X(2) \oplus \dots )
\]
  The vertex $p$ of
$C_X$ is the (closed) subscheme defined by the augmentation ideal,
$ker\:{\epsilon}$, where
\[
  \epsilon: ({\cal O}_X \oplus
{\cal O}_X(1) \oplus {\cal O}_X(2) \oplus \dots) \rightarrow {\cal
O}_X
\]
  $\Spec ({\epsilon})$ thus defines a map $X \rightarrow C_X$ whose image is
$p$.\end{r}

For future use we recall that the cone with vertex removed
(\'{e}point\'{e})
\[C_X - p \simeq \V({\cal O}_X(-1))\simeq \Spec (\oplus {\cal O}_X(n))\]
 where $\V$ denotes
the operation of taking the vector bundle associated to a locally free
sheaf.  There is an action of $\G_m$ on $C_X$ (and on $C_X-p$) with integral
weights coming from the action on each tensor power ${\cal O}(n)$.
For details see [{\bf {EGA II}} \S 8.4 - 8.6].

Next we define a deformation functor $\Def_{C_X}$
as the functor which assigns to each object $A$ in $\cal C$ the set of
deformations ${\cal O}_{C_X,p}(A)$ of the $k$-algebra ${\cal
O}_{C_X,p}$ which is the local ring of $C_X$ at $p$.

There is a natural morphism  $h$ from the Hilbert functor to this deformation
functor obtained by assigning to a scheme $Y$ as above the local ring
of the vertex of the projective cone on
$Y$.  This is naturally a deformation of the local ring of the vertex
of the projective cone on $X$.  The morphism thus consists of ``forgetting the
embedding in $\Pj^n$
or, in terms of the underlying rings, forgetting the gradings.

Comparison Theorems between Hilbert functors and Deformation functors
go back (at least) to Schlessinger \cite{Sch:71} and can be found in
\cite{P:74}, \cite{K:79} and \cite{Wa:92}. We will use a version due
to Kleppe which employs Andr\'{e}-Quillen cohomology to avoid unnecessary
smoothness conditions.

Since $X$ is given as the zero-scheme of a section of a vector bundle we want
to describe its Hilbert scheme in the same way.  So define functors
$F_{s_0}$ from the category ${\cal C}$ to {\em Sets} which, for any given
section
$s_0 \in H^0(Z,E)$, assign to an object $A$ of ${\cal C}$ the set of zero
schemes in $Z \times \Spec A$  of sections $s_A \in H^0(Z \times\Spec A,
E \times \Spec A)$ which reduce to $s_0$ over the closed point $k$ under
the  map $\Spec k \rightarrow \Spec A$.

If we further assume that codim$X$=rank$E$ we have
that $X$ is (at least) a local complete intersection.  Moreover, such zero
schemes $s_A$
are flat over
$\Spec A$. By [{\bf EGA IV}
\S  19.3.8]  the zero scheme defined by any of
the $s_A$ is a
local complete intersection as well.  Thus we obtain elements of ${\Hilb}_X(A)$
as sections of
vector bundles $E \times \Spec A$ and a morphism of functors $F_{s_0}
\stackrel{f}{\rightarrow} \Hilb_X$.

We now study this morphism of functors.

\begin{s}{Proposition}
\label{surj}
 Let $Z$, $E$ and $A$ be as above and let $X$ be
the zero set of a section $s_0$ with codim $X$=rank$E$, write ${\cal
I}_X \subset {\cal  O}_Z$ for the ideal sheaf of $X$ and suppose that
\[H^1(Z,E \otimes{\cal I}_X)=0\]
then the above morphism of functors
$f: F \rightarrow \Hilb$ is surjective on tangent spaces.\end{s}

\pf.   The tangent space to $H$ at $h_0$, is
$\Hilb_{s_0}(k[{\epsilon}])$ and it
is standard that
\[T(H,h_0)=Hom_{{\cal O}_Z}({\cal I}_X, {\cal O}_X)=Hom_{{\cal
O}_X}({\cal I}_X/{\cal I}^2_X, {\cal O}_X)\]
(for the first isomorphism see \cite {G:61}, for the second see
\cite{H:77} Ch. II, \S 8).  We also have
\[Hom_{{\cal
O}_X}({\cal I}_X/{\cal I}^2_X, {\cal O}_X)\simeq H^0(X,N_{X/Z}) \simeq
H^0(X,E\!\mid_X)\]
where $N$ is isomorphic to the tangent space
to $\Hilb$.  But the tangent space to
$F:=F(k[{\epsilon}])$ is the vector space of sections of
$E \times \Spec(k[{\epsilon}]) \rightarrow Z \times
\Spec(k[{\epsilon}]) $ which reduce to $s_0$
over $k$ and this is just $H^0(Z,E)$.
Consider the standard short exact sequence
\[
0 \rightarrow {\cal I}_X \rightarrow {\cal O}_Z \rightarrow {\cal O}_X
\rightarrow 0
\]
tensor it with $E$ and take cohomology to get
\[0 \rightarrow H^0(Z,E\otimes {\cal I}_X) \rightarrow H^0(Z,E)
\rightarrow H^0(X,E\!\mid_X)\]
\[
\rightarrow H^1(Z,E\otimes {\cal I}_X) \rightarrow \dots \hfill
\]
Hence if $H^1(Z,E\otimes{\cal I}_X)=0$ we have the desired surjectivity.
 \fp

So far we have developed portions of a triangle of functors

\commtriang{$F_{s_0}$}{$\Hilb_{s_0}$}{$\Def_{C_X}$}{$f$}{?}{h}

To fill in the morphism marked with ``?'' note that the local ring of the
vertex of the cone of the zero scheme of $s_A$ is naturally a
deformation of the local ring of the vertex of the cone of the zero
scheme of $s_0$.  It is clear that the triangle commutes because $g$
is just the composition of $f$ and $h$.  We now
invoke the comparison Theorem relating the Hilbert functor and the deformation
functor

\begin{s}{Theorem}  Suppose that $V$ is projectively normal and that
$T^1(V)$,the tangent space to $\Def_{C_V}$, is
negatively graded (in a sense to be made precise below), then the
natural morphism of functors
\[h: \Hilb_{V} \rightarrow \Def_{C_V}\]
is smooth.
\end{s}

\pf.  See \cite{K:79}  \fp

As we will see in the case of the zero schemes of sections of
certain vector bundles, this result together with Proposition \re{surj} will
enable us to establish the
existence of families of the desired type.  The codim$X$=rank$E$
condition that we have been placing on our zero schemes of sections
of $E$ ensures projective normality (assuming that the varieties are normal in
the first place)
by a straightforward adaptation of the proof for complete
intersections in projective space.  It is the grading condition in the
comparison theorem which
is the more difficult of the two to verify and we will reduce it to a
condition on vanishing of cohomology.

Recall that $\Def_{C_X}$  is the deformation functor of the vertex on the
projective cone $C_X$ over $X$ and that there is an  action by $\G_m (k)$ on
$C_X$
with weights ranging over the integers.
$T^1(X):={\Def}_{C_X}(k[{\epsilon}])$ is a vector space and the action
of $\G_m$ on on $C_X$ becomes an action on $T^1(X)$, so that we get a
decomposition $T^1(X)=\oplus_{\nu = -\infty}^{\infty}T^1(\nu )$ and
$T^1(X)$ becomes a graded vector space (as $p$ is an isolated
singularity in $C_X$, in fact $T^1(\nu)=0$ for $\nu \gg 0$ and for $\nu \ll 0$
but we will not need this).  We say that $T^1(X)$ is
{\em negatively graded} if $T^1(X)(\nu )=0$ for $\nu > 0$.  To demonstrate that
this condition obtains
we show that $\Hilb (k[{\epsilon}])$ too has a grading,  which is negative in
this
sense and that the map $Th$ is surjective and respects both gradings.

\begin{s}{Proposition} Let $X$ be the zero scheme of the section $s_0$ of $E
\rightarrow Z$ and assume that $X$ projectvely
normal (e.g. codim $X$ = rank $E$).
Suppose the conclusion of Proposition \re{surj} holds, that is, the morphism of
functors $F_{s_0} \rightarrow
\Hilb_X$ is surjective on tangent spaces. Let $g: F_{s_0}
\rightarrow Def_{C_X}$ be the natural morphism taking an element of
$F_{s_0}(A)$ to an element of $Def_{C_x}(A)$.  Denote the map on
tangent spaces by $H^0(Z,E) \stackrel{Tg}{\rightarrow} T^1(X)$.  If
this map is surjective then $T^1(X)$ is negatively graded and the morphism of
functors $h$ is smooth.
\end{s}

\pf.  We have a triangle of tangent maps

\commtriang{$H^0(Z,E)$}{$\Hilb(k[{\epsilon}])$}{$\Def{C_X}(k[{\epsilon}])$}{$Tf$}{$Tg$}{$Th$}

which commutes because the triangle of underlying maps does and we are
assuming that $Tf$ is onto.  If we now further assume that
$H^0(Z,E) \rightarrow T^1(X)$ is onto then $\Hilb(k[{\epsilon}])
\rightarrow T^1(X)$ must also be onto as well.

$\Hilb(k[{\epsilon}])\simeq H^0(X,N_{X/Z})$ as noted above and
$H^0(X,N_{X/Z}) \simeq H^0(C_X,N_{C_X})$ where $C_X$ is the projective cone
over $X$.  By projective normality the singularity at the vertex $p$ of
$C_X$ has depth=2 so global sections of $C_X-p$ extend to global sections
of $C_X$.  Since $X$ is a local complete intersection, this is true of
all coherent sheaves on $C_X$ by local duality, and so, in particular
\[H^0(C_X,N_{C_X}) \simeq H^0({C_X}-p,N_{C_X}) \simeq H^0(C_X-p,N_{C_X-p})\]

But
\[C_X-p \simeq \V({\cal O}_X(-1)) - {\em zero section}\]
 so there is a natural affine
map $\pi:C_X-p \rightarrow X$ and $N_{C_C-p} \simeq \pi^*N_X$ and we can
compute
\[H^0(C_X-p,N_{C_X-p}) \simeq H^0(C_X-p,\pi^*N_X) \simeq
H^0(X,\pi_*\pi^*N_X) \simeq \oplus_{\nu = -\infty}^0 H^0(X,N_X(\nu))\]
where $N_X(\nu):= N_X \otimes {\cal O}(\nu)$ as usual.
The second isomorphism comes from the fact that $\pi$ is affine and
the conclusion that $H^0(X,N_X(\nu))=0$ for $\nu > 0$ comes from the fact that
the sections must extend over all of $C_X$ and hence cannot have any
poles at $p$.

The grading thus produced on $\Hilb (k[{\epsilon}])$ arises from the action of
$\G_m$ on $C_X$ and the morphism $\Hilb \rightarrow \Def$ is contructed
precisely by passing to
the projective cone to go from the varieties in projective space to
abstract varieties.  Hence the action of $\G_m$ on $T^1(X)$ and
$\Hilb(k[{\epsilon}])$ is compatible with this morphism.  The result is that
$T^1(X)$ must also be negatively graded and hence we can apply the
comparison Theorem to get that this morphism is smooth.  \fp

Finally it is not difficult to reduce the requirement ``$H^0(Z,E)
\rightarrow T^1(X)$ surjective'' to a statement of vanishing of
cohomology.

The tangent sheaf to a variety $V$ is defined by $\Theta_V :=
Hom_{{\cal O}_V}(\Omega_V, {\cal O}_V)$, where $\Omega_V$ is the sheaf
of differentials of $V$.  Since $C_X \subset \Aff^{n+1}$ and $C_X -p$
is smooth we have that
\[
0 \rightarrow \Theta_{C_X}\!\mid_{C_X-p}\rightarrow \Theta_{\Aff
^{n+1}}\!\mid_{C_X-p} \rightarrow {\cal N}_{C_X}\!\mid_{C_X-p}
\rightarrow 0
\]
is exact.  Note that $C-p$ is not affine so there is a longer exact
cohomology sequence
\[
0 \rightarrow  H^0(C_X-p, \Theta_{C_X}\!\mid_{C_X-p})
\rightarrow  H^0(C_X-p,\Theta_{\Aff^{n+1}\!\mid_{C_X-p}})
\rightarrow\]
\[H^0(C_X-p,{\cal N}_{C_X}\!\mid_{C_X-P})  \rightarrow  H^1(C_X-p,
\Theta_{C_X}\!\mid_{C_X-p})  \rightarrow H^1(C_X-p,
\Theta_{\Aff^{n+1}\!\mid_{C_X-p}})  \rightarrow
\]

The theory of the ``cotangent complex'' \cite{Li-Sch 67} gives us the following
exact sequence
\[
0 \rightarrow H^0(C,\Theta_{C_X}) \rightarrow
H^0(C,\Theta_{\Aff^{n+1}}\!\mid_{C_X}) \rightarrow
H^0(C_X,{\cal
N}_{C_X})\]
\[
 \rightarrow T^1(X) \rightarrow H^1(X,\Theta_Z\!\mid_X) \rightarrow
\dots \hfill
\]

\noindent Once again since $X$ is projectively normal, sections over $C_X-p$
extend to $C_X$ so that
\[H^0(C_X,\Theta_{C_X}) \simeq
H^0(C_X-p,\Theta_{C_X}\!\mid_{C_X-p})
\]

\[
H^0(C_X,\Theta_{\Aff^{n+1}}\!\mid_{C_X}) \simeq
H^0(C_X-p,\Theta_{\Aff^{n+1}}\!\mid_{{C_X}-p})\]
and
\[H^0(C_X,{\cal N}_{C_X}) \simeq H^0(C_x-p,{\cal N}_{C_X-p})\]

Putting all of this together we see that if
\[
 H^1(C_x-p,\Theta_{\Aff^{n+1}}\!\mid_{C_X-p})\simeq
H^1(X,\Theta_Z\!\mid_X)=\varnothing
\]
then
\[
 T^1(X) \hookrightarrow H^1(C_X-p,\Theta_{C_X}\!\mid_{C_X-p})=0
\]

To summarize,

\begin{s}{Proposition} For $X$ the zero scheme of a section $S$ bundle $E$ over
a non-singular
projective variety $Z$ with codim $X$ = rank $E$ where
\[H^1(Z, E
\otimes {\cal I}_X)=H^1(X,\Theta_Z\!\mid_X)=0\]
and $X$ is normal we have the commutative triangle

\commtriang{$F_{s_0}$}{$\Hilb_X$}{$\Def_{C_X}$}{}{}{}

where all of the arrows are smooth. \end{s}

\pf.  We only need to show that $F_{s_0} \stackrel{f}{\rightarrow} \Hilb_X$ is
smooth since then the third side of the triangle will be smooth by
\cite{Sch:67}. We know that this arrow is surjective on tangent spaces
hence we only need to show that $F_{s_0}$ is less obstructed than
$\Hilb_X$.  In other words, if $B \stackrel{\phi}{\rightarrow} A$ is a
surjection in
$\cal C$ with $ker(\phi)^2=0$ and $({\bf m}_B)ker(\phi)=0$ (so $ker(\phi)$ is a
$k$-vector space) and if for $\xi_0 \in \Hilb_X(A)$ which is
in the image of $f$ so, $\xi_0=f(\zeta_0)$, $\zeta_0 \in F_{s_0}(A)$,
there is a $\xi \in \Hilb_X(B)$ such that $\Hilb_X(\phi)(\xi)=\xi_0$
then there is a $\zeta \in F_{s_0}(B)$ such that
$F_{s_0}(\phi)(\zeta)=\zeta_0$.  The obstruction to the existence of
such $\zeta$ lies in $Ext^2(L_{X/Z},{\cal O}_X)$ where $L_{X/Z}$ is
the cotangent complex of $X$.  This is a two term complex
\[
0 \rightarrow {\cal I}_X/({\cal I}_X)^2 \rightarrow i^*(\Omega^1_Z)
\rightarrow 0
\]
since $X$ is a local complete intersection \cite{I:71} Ch. III, \S 3.2 and the
$Ext^2$ terms vanish
since $H^1(X,\Theta_Z\!\mid_X)=0$. \fp \bl

Since we know that the formal scheme prorepresenting $\Hilb_X$ is
algebraizable \cite{G:61}, we now have the same conclusion for $\Def_{C_X}$.
All
infinitesimal deformations of $X$ come from small deformations and all
small deformations lie in $Z$.  Indeed, by versality
$H^0(Z,E)$ is a complete deformation space and the family
\[
Y:=\{(z,s)\mid z \in Z, s \in H^0(Z,E), s(z)=0\} \subset Z \times
H^0(Z,E) \rightarrow H^0(Z,E)
\]
is a universal family both in the sense of the deformation functor and
the Hilbert functor.  These results are valid for $X$ defined over any field.

\bitt{Application of the theory}

We now apply this theory to our problem.  If a normal $k$-variety $V$ ($k$ a
finitely generated extension of $\Q$) is
given as the zero scheme of a section $s\in H^0(Z,E)$ (everything
defined over $k$) and if we have
\begin{itemize}
\item codim$V$=rank$E$
\item $H^1(Z,E\otimes {\cal I}_V)=H^1(V,\Theta_V)=0$
\end{itemize}
then we know that the $k$-vector space $H^0(Z,E)$ parametrizes a
complete family of deformations of $V$.  There is a map
\[
H^0(Z,E) \rightarrow H^0(Z,E)_{/{\Q}} \otimes k
\]
giving a $\Q$ structure to $H^0(Z,E)$.  This induces a morphism of
functors
$F^k_{s} \rightarrow F^{\Q}_{s}$.

We observe that this morphism is
smooth since first, if $A' \rightarrow A$ is a surjection in $\cal C$ then
\[
F^k(A') \rightarrow F^k(A) \times_{F^{\Q}(A)} F^{\Q}(A')
\]
must also be onto.  To see this note that the arrow $\beta$ in
\[
\begin{array}{ccc}
 F^k(A) \times_{F^{\Q}(A)} F^{\Q}(A') & \stackrel{{\beta}'}{\rightarrow} &
F^{\Q}(A')\\
\downarrow &     & \downarrow\\
F^k(A) & \stackrel{\beta}{\rightarrow} & F^{\Q}(A)
\end{array}
\]
is onto and hence ${\beta}'$ is also and that $F^k(A') \rightarrow
F^{\Q}(A')$ is onto as well.  Secondly the morphism also induces a
bijection on tangent spaces since $H^0(Z,E)_k$ and $H^0(Z,E)_{/{\Q}}$
are vector spaces of the same dimension over fields of the same cardinality.

Thus the projection $H^0(Z,E) \rightarrow H^0(Z,E)_{/{\Q}}$ induces a map of
deformation spaces between the functors $\Def^k$ and ${\Def}^{\Q}$ and
$\Hilb_X^k$ and $\Hilb_X^{\Q}$.

We have a universal family $Y_{/{\Q}} \rightarrow H^0(Z,E)_{/{\Q}}$
and the above ensures that we recover the corresponding universal
family over $k$ by extension of scalars
\[
\begin{array}{ccc}
Y_k & \rightarrow & Y_{/{\Q}}\\
\downarrow &   & \downarrow \\
H^0(Z,E)_k & \rightarrow & H^0(Z,E)_{/{\Q}}
\end{array}
\]
But there is a unique extension of scalars from $\Q$ to $\C$ and this
gives  a diagram of the
sort described at the beginning of this \S

\commrect{$Y_{/{\C}}\simeq Y_{/{\Q}} \times
\C$}{$Y_{/{\Q}}$}{$H^0(Z,E)_{/{\C}}\simeq
H^0(Z,E)_{/{\Q}}\times
\C$}{$H^0(Z,E)_{/{\Q}}$}{$f$}{$f_{/{\Q}}$}{$\beta_Y$}{$\beta$}

Since we can now consider $V$, $E$, $Z$, and $s$ as objects defined
over the complex numbers it is clear that all of the conjugate
varieties $V_{\sigma}$ can be obtained by conjugating the section $s$, that is,
$V_{\sigma}$ is the zero scheme of $s_{\sigma}:=\sigma (s)$.  Hence all
of the conjugates $V$ are in the family $Y$ if $V$ is.  Thus the
conditions of Theorem \re{main} are satisfied and the independence of
the topological type from variation with $\sigma$ is ensured.

It only remains to spell out some types of varieties which are defined as the
zero sections of vector bundles satisfying our conditions.

If we assume that $Z$ is homogeneous, that is
\[Z \simeq G/H\]
where $G$ is a simple, simply connected, split algebraic group over
$\Q$, and $H$ a parabolic subgroup, then by algebraic versions of Bott's
vanishing Theorems \cite{Dm:76} we get that $H^1(Z, \Theta_Z)=0$ for all $i >
0$.
If we further suppose that $E\simeq \oplus_{j=1}^r{\cal O}_Z(d_j)$
with $d_j > 0$ and $X$ defined by a section $s$ such that
\begin{itemize}
\item codim $X$=$r$
\item $X$ is not a $K$-3 surface
\item dim $Z \geq 3$
\end{itemize}
then a reworking of calculations of Borcea using an algebraic version of
the Kodaira-Nakano-Akizuki vanishing theorem shows that
$H^1(Z,E\otimes{\cal I}_X)=0$.  For details see \cite{Re:94}

If we now assume that $X$ is normal and defined over $k$ that is, the section
$s$ satisfies
\[
s \in H^0(Z,E)_{/{\Q}} \otimes k
\]
\noindent then we have a class of $k$ varieties whose topology
is independent of the embedding of $k$. This class includes complete
intersections in projective space.  A somewhat different proof is available for
the special case of complete intersections in projective space which does not
require normality \cite{Re:94}.

We note finally that while the case of
$K$-3 surfaces must be treated separately, the result is the same
since all $K$-3's are homeomorphic (indeed diffeomorphic) and
conjugates of $K$-3 varieties are $K$-3.

\bit{Comments and Open Questions}

The above discussion can be used to shed a bit of light on the nature of the
Serre-Abelson examples.  Both authors construct their varieties as
quotients using finite group actions.  In both cases the action is
varied under conjugation.  The difference lies in the part of the
homotopy type which is affected by conjugation.

In Serre, the variety acted upon is a product of a diagonal
hypersurface by an abelian variety and the group action on the abelian
variety makes the $\pi_1$ of the variety into a module  in
demonstrably different ways under conjugation.  In Abelson's example
the group acts on complete intersection which is constructed via a
representation of the group and all of this varies under conjugation.
Special choices of the group allow one to demonstrate variation in the
Postnikov tower.

It is not difficult to see, in both cases that the field of moduli of the
varieties thus constructed is smaller than their field of definition.  The
group action creates some ``symmetries'' which produce this result.

This suggests the following vague question:

A)  Is it possible to produce examples of non-homeomorphic conjugate
varieties without using (finite) group actions?

This question can be made more specific in a number of ways.

Because of the use of group actions the examples of Serre and Abelson
are rather rigid.  A small deformation of one of their examples no
longer maintains the structure required to compare it with its
conjugates.  One approach might therefore be to ask,

A1) Is it possible to construct examples of non-homemorphic conjugate
varieties which are stable under small deformations?  This seems unlikely.

Another way to cut out finite group actions is to ask for
simply-connected examples,

A2)  Are there examples of simply connected non-homemorphic conjugate
varieties?

A useful source of new examples may be provided by Shimura varieties.

Note: It may be useful to employ Kollar's
notions of ''essentially large'' fundamental groups here instead.

Finally, along these lines one has the fundamental question,

A4) Are the simply connected covering spaces of conjugate algebraic varieties
analytically isomorphic?

The criterion developed in this paper does not provide an indication
of the minimum '``necessary'' conditions under which the topology of varieties
remains stable under conjugation.  Neither does it give any indication of the
``part'' (if any) of the topological type which are conjugations invariant
(over and above the \'{e}tale homotopy type which is clearly invariant).  A
Theorem of Deligne \cite{D
e:87} shows that the nilpotent completion of the fundamental group of an
algebraic variety is algebraically determined.   One is led to ask,

B)  Is the entire rational homotopy type a conjugation invariant?

One might also pose the following question which seems to lie somewhere between
A) and B),

C)  Is simple connectivity a conjugation invariant?

\end{enumerate}
\end{document}